\documentclass{article}
\usepackage{eurosym}
\usepackage{amssymb, amsmath}

\begin{document}

\begin{center}
{\large \ Properties and applications of some special integer number
sequences }%
\begin{equation*}
\end{equation*}

\textbf{Cristina Flaut, Diana Savin and Geanina Zaharia}%
\begin{equation*}
\end{equation*}
\end{center}

\textbf{Abstract.} {\small In this paper, we provide properties and
applications of some special integer sequences. We generalize and give some
properties of Pisano period. Moreover, we provide a new application in
Cryptography and applications of some quaternion elements. } 
\begin{equation*}
\end{equation*}


\medskip

\begin{equation*}
\end{equation*}

\begin{center}
\bigskip
\end{center}

There have been numerous papers devoted to the study of the properties and
applications of particular integer sequences, the most studied of them being
Fibonacci sequence ([FP; 09], [Ha; 12], \textbf{[}Ho; 63\textbf{],} [St;
06], [St; 07], etc.). In this paper, we study properties and provide some
applications of these number sequences in a more general case. 

Let $a_{1},...,a_{k}$ be arbitrary integers, $a_{k}\neq 0$. We consider the
general $k-$terms recurrence, $n,k\in \mathbb{N},k\geq 2,n\geq k,$%
\begin{equation}
d_{n}\text{=}a_{1}d_{n-1}\text{+}a_{2}d_{n-2}\text{+}...\text{+}%
a_{k}d_{n-k},d_{0}\text{=}d_{1}\text{=...=}d_{k-2}\text{=}0,d_{k-1}\text{=}1,
\tag{1.1.}
\end{equation}%
where $a_{k}\neq 0$ are given integers and its associated matrix $D_{k}\in 
\mathcal{M}_{k}\left( \mathbb{R}\right) ,$ 
\begin{equation}
D_{k}=\left( 
\begin{array}{ccccc}
a_{1} & a_{2} & a_{3} & ... & a_{k} \\ 
1 & 0 & 0 & ... & 0 \\ 
0 & 1 & 0 & ... & 0 \\ 
... & ... & ... & ... & ... \\ 
0 & 0 & ... & 1 & 0%
\end{array}%
\right) ,  \tag{1.2.}
\end{equation}%
see [Jo].

\text{For }$n\in \mathbb{Z},n\geq 1,$\text{ we have that }

\begin{equation}
D_{k}^{n}\text{=}\left( 
\begin{array}{ccccc}
d_{n\text{+}k\text{-}1} & \overset{k\text{-}1}{\underset{i\text{=}1}{\sum }}%
a_{i\text{+}1}d_{n\text{+}k\text{-}i\text{-}1} & \overset{k-2}{\underset{i=1}%
{\sum }}a_{i\text{+}2}d_{n\text{+}k\text{-}i\text{-}1} & \text{...} & 
a_{k}d_{n\text{+}k\text{-}2} \\ 
d_{n\text{+}k\text{-}2} & \overset{k\text{-}1}{\underset{i\text{=}1}{\sum }}%
a_{i\text{+}1}d_{n\text{+}k\text{-}i\text{-}2} & \overset{k-2}{\underset{i=1}%
{\sum }}a_{i\text{+}2}d_{n\text{+}k\text{-}i\text{-}2} & \text{...} & 
a_{k}d_{n\text{+}k\text{-}3} \\ 
d_{n\text{+}k\text{-}3} & \overset{k\text{-}1}{\underset{i\text{=}1}{\sum }}%
a_{i\text{+}1}d_{n\text{+}k\text{-}i\text{-}3} & \overset{k-2}{\underset{i=1}%
{\sum }}a_{i\text{+}2}d_{n\text{+}k\text{-}i\text{-}3} & \text{...} & 
a_{k}d_{n\text{+}k\text{-}4} \\ 
\text{...} & \text{...} & \text{...} & \text{...} & \text{...} \\ 
d_{n} & \overset{k\text{-}1}{\underset{i\text{=}1}{\sum }}a_{i\text{+}1}d_{n%
\text{+}k\text{-}i\text{-}k} & \overset{k-2}{\underset{i=1}{\sum }}a_{i\text{%
+}2}d_{n\text{-}i} & \text{...} & a_{k}d_{n\text{-}1}%
\end{array}%
\right) \text{,}  \tag{1.3.}
\end{equation}%
see [Fl; 19], Proposition 2.1.

In [Me; 99], the author proved some new properties of Fibonacci and Pell
numbers. In the following, we will generalize some of these results.

If we consider the matrix $Y_{i}=\left( d_{i+k-1},...,d_{i+1},d_{i}\right) $%
, for $i=0$, relation $\left( 1.1\right) $ can be written under the form 
\begin{equation}
Y_{1}^{t}=D_{k}Y_{0}^{t}.  \tag{1.4.}
\end{equation}

\textbf{Proposition 1.1.} \textit{With the above notations, we have the
following relations:}

\begin{equation}
Y_{n}^{t}=D_{k}Y_{n-1}^{t}  \tag{1.5.}
\end{equation}%
\begin{equation}
Y_{n}^{t}=D_{k}^{n}Y_{0}^{t}.  \tag{1.6.}
\end{equation}%
\begin{equation}
Y_{n+r}^{t}=D_{k}^{n}Y_{r}^{t}.  \tag{1.7.}
\end{equation}

\textbf{Proof}.

Relation $\left( 1.5\right) $ is obviously.

To prove relation $\left( 1.6\right) ,$ we remark that for $n=1$, we have
relation $\left( 1.4\right) .$ Supposing relation true for $n-1,$ we prove
it for $n.$ From relation $\left( 1.5\right) ,$ we have that $%
Y_{n}^{t}=D_{k}Y_{n-1}^{t}=D_{k}^{n}Y_{0}^{t}$.

To prove relation $\left( 1.7\right) $, we have $%
Y_{n+r}^{t}=D_{k}^{n+r}Y_{0}^{t}=D_{k}^{n}D_{k}^{r}Y_{0}^{t}=D_{k}^{n}Y_{r}^{t}.\Box \medskip 
$\medskip

\textbf{Proposition 1.2}. \textit{With the above notations, for difference
equation }$\left( 1.1\right) $\textit{, the following relation are true:}

i)%
\begin{equation*}
\left\vert 
\begin{array}{ccccccc}
d_{n+k-1} & d_{n+k} & d_{n+k+1} & \text{...} & d_{n+2k-4} & d_{n+2k-3} & 1
\\ 
d_{n+k-2} & d_{n+k-1} & d_{n+k} & \text{...} & d_{n+2k-5} & d_{n+2k-4} & 0
\\ 
d_{n+k-3} & d_{n+k-2} & d_{n+k-1} & \text{...} & d_{n+2k-6} & d_{n+2k-5} & 0
\\ 
... & ... & ... & \text{...} & ... & ... & \text{...} \\ 
d_{n+2} & d_{n+3} & d_{n+4} & \text{...} & d_{n+k-2} & d_{n+k-1} & 0 \\ 
d_{n+1} & d_{n+2} & d_{n+3} & \text{...} & d_{n+k-3} & d_{n+k-2} & 0 \\ 
d_{n} & d_{n+1} & d_{n+2} & \text{...} & d_{n+k-4} & d_{n+k-3} & 0%
\end{array}%
\right\vert \text{=}\left( \text{-}1\right) ^{n\left( k+1\right)
}a_{k}^{n}d_{-n}.
\end{equation*}

ii)%
\begin{equation*}
\left\vert 
\begin{array}{cccccc}
d_{n+k-1} & d_{n+k} & d_{n+k+1} & ... & d_{n+2k-4} & d_{n+2k-3} \\ 
d_{n+k-2} & d_{n+k-1} & d_{n+k} & ... & d_{n+2k-5} & d_{n+2k-4} \\ 
d_{n+k-3} & d_{n+k-2} & d_{n+k-1} & ... & d_{n+2k-6} & d_{n+2k-5} \\ 
... & ... & ... & ... & ... & ... \\ 
d_{n+2} & d_{n+3} & d_{n+4} & ... & d_{n+k-1} & d_{n+k} \\ 
d_{n+1} & d_{n+2} & d_{n+3} & ... & d_{n+k-2} & d_{n+k-1}%
\end{array}%
\right\vert \text{=}\left( \text{-}1\right) ^{n\left( k+1\right) }a_{k}^{n}%
\text{.}
\end{equation*}

iii) 
\begin{equation*}
d_{m\text{+}k\text{-}1}d_{n\text{+}k\text{-}1}\text{+}d_{m\text{+}k\text{-}2}%
\overset{k-1}{\underset{i=1}{\sum }}a_{i\text{+}1}d_{n\text{+}k\text{-}i%
\text{-}1}\text{+}d_{m+k-3}\overset{k-2}{\underset{i=1}{\sum }}a_{i\text{+}%
2}d_{n\text{+}k\text{-}i\text{-}1}\text{+...+}a_{k}d_{m}d_{n\text{+}k\text{-}%
2}\text{=}d_{m\text{+}n\text{+}k\text{-}1}\text{.}
\end{equation*}

\textbf{Proof.} 1) Indeed, since $\det D_{k}=\left( -1\right) ^{k+1}a_{k},$%
we obtain that \newline
\begin{equation*}
\left\vert 
\begin{array}{ccccccc}
d_{n+k-1} & d_{n+k} & d_{n+k+1} & \text{...} & d_{n+2k-4} & d_{n+2k-3} & 1
\\ 
d_{n+k-2} & d_{n+k-1} & d_{n+k} & \text{...} & d_{n+2k-5} & d_{n+2k-4} & 0
\\ 
d_{n+k-3} & d_{n+k-2} & d_{n+k-1} & \text{...} & d_{n+2k-6} & d_{n+2k-5} & 0
\\ 
... & ... & ... & \text{...} & ... & ... & \text{...} \\ 
d_{n+2} & d_{n+3} & d_{n+4} & \text{...} & d_{n+k-1} & d_{n+k} & 0 \\ 
d_{n+1} & d_{n+2} & d_{n+3} & \text{...} & d_{n+k-2} & d_{n+k-1} & 0 \\ 
d_{n} & d_{n+1} & d_{n+2} & \text{...} & d_{n+k-3} & d_{n+k-2} & 0%
\end{array}%
\right\vert \text{=}\newline
\end{equation*}%
\begin{equation*}
\text{=}\left\vert D_{k}^{n}\left( 
\begin{array}{c}
d_{k-1} \\ 
d_{k-2} \\ 
d_{k-3} \\ 
\text{...} \\ 
d_{2} \\ 
d_{1} \\ 
d_{0}%
\end{array}%
\right) D_{k}^{n}\left( 
\begin{array}{c}
d_{k} \\ 
d_{k-1} \\ 
d_{k-2} \\ 
\text{...} \\ 
d_{3} \\ 
d_{2} \\ 
d_{1}%
\end{array}%
\right) D_{k}^{n}\left( 
\begin{array}{c}
d_{k+1} \\ 
d_{k} \\ 
d_{k-1} \\ 
\text{...} \\ 
d_{4} \\ 
d_{3} \\ 
d_{2}%
\end{array}%
\right) \text{...}D_{k}^{n}D_{k}^{-n}\left( 
\begin{array}{c}
d_{\text{-}n\text{+}k\text{-}1} \\ 
d_{\text{-}n\text{+}k\text{-}2} \\ 
d_{\text{-}n\text{+}k\text{-}3} \\ 
\text{...} \\ 
d_{\text{-}n\text{+}2} \\ 
d_{\text{-}n\text{+}1} \\ 
d_{\text{-}n}%
\end{array}%
\right) \right\vert \text{=}\newline
\end{equation*}%
\begin{equation*}
\text{=}\det \left( D_{k}^{n}\right) \left\vert 
\begin{array}{ccccccc}
d_{k-1} & d_{k} & d_{k+1} & \text{...} & d_{2k-4} & d_{2k-3} & d_{-n+k-1} \\ 
d_{k-2} & d_{k-1} & d_{k} & \text{...} & d_{2k-5} & d_{2k-4} & d_{-n+k-2} \\ 
d_{k-3} & d_{k-2} & d_{k-1} & \text{...} & d_{2k-6} & d_{2k-5} & d_{-n+k-3}
\\ 
... & ... & ... & \text{...} & ... & ... & ... \\ 
d_{2} & d_{3} & d_{4} & \text{...} & d_{k-1} & d_{k} & d_{-n+2} \\ 
d_{1} & d_{2} & d_{3} & \text{...} & d_{k-2} & d_{k-1} & d_{-n+1} \\ 
d_{0} & d_{1} & d_{2} & \text{...} & d_{k-3} & d_{k-2} & d_{-n}%
\end{array}%
\right\vert \text{=}\left( \text{-}1\right) ^{n\left( k+1\right)
}a_{k}^{n}d_{-n}.
\end{equation*}

2) Indeed, using the ideas developed above, obtain\newline
\begin{equation*}
\left\vert 
\begin{array}{cccccc}
d_{n+k-1} & d_{n+k} & d_{n+k+1} & ... & d_{n+2k-4} & d_{n+2k-3} \\ 
d_{n+k-2} & d_{n+k-1} & d_{n+k} & ... & d_{n+2k-5} & d_{n+2k-4} \\ 
d_{n+k-3} & d_{n+k-2} & d_{n+k-1} & ... & d_{n+2k-6} & d_{n+2k-5} \\ 
... & ... & ... & ... & ... & ... \\ 
d_{n+2} & d_{n+3} & d_{n+4} & ... & d_{n+k-1} & d_{n+k} \\ 
d_{n+1} & d_{n+2} & d_{n+3} & ... & d_{n+k-2} & d_{n+k-1}%
\end{array}%
\right\vert =
\end{equation*}%
\begin{equation*}
\text{=}\left\vert D_{k}^{n}\left( 
\begin{array}{c}
d_{k-1} \\ 
d_{k-2} \\ 
d_{k-3} \\ 
\text{...} \\ 
d_{2} \\ 
d_{1}%
\end{array}%
\right) D_{k}^{n}\left( 
\begin{array}{c}
d_{k} \\ 
d_{k-1} \\ 
d_{k-2} \\ 
\text{...} \\ 
d_{3} \\ 
d_{2}%
\end{array}%
\right) D_{k}^{n}\left( 
\begin{array}{c}
d_{k+1} \\ 
d_{k} \\ 
d_{k-1} \\ 
\text{...} \\ 
d_{4} \\ 
d_{3}%
\end{array}%
\right) \text{...}D_{k}^{n}\left( 
\begin{array}{c}
d_{2k-4} \\ 
d_{2k-5} \\ 
d_{2k-6} \\ 
\text{...} \\ 
d_{k-1} \\ 
d_{k-2}%
\end{array}%
\right) D_{k}^{n}\left( 
\begin{array}{c}
d_{2k-3} \\ 
d_{2k-4} \\ 
d_{2k-5} \\ 
\text{...} \\ 
d_{k} \\ 
d_{k-1}%
\end{array}%
\right) \right\vert \text{=}\newline
\end{equation*}%
\begin{equation*}
\text{=}\det \left( D_{k}^{n}\right) \left\vert 
\begin{array}{cccccc}
d_{k-1} & d_{k} & d_{k+1} & ... & d_{2k-4} & d_{2k-3} \\ 
d_{k-2} & d_{k-1} & d_{k} & ... & d_{2k-5} & d_{2k-4} \\ 
d_{k-3} & d_{k-2} & d_{k-1} & ... & d_{2k-6} & d_{2k-5} \\ 
... & ... & ... & ... & ... & ... \\ 
d_{2} & d_{3} & d_{4} & ... & d_{k-1} & d_{k} \\ 
d_{1} & d_{2} & d_{3} & ... & d_{k-2} & d_{k-1}%
\end{array}%
\right\vert \text{=}\left( \text{-}1\right) ^{n\left( k+1\right) }a_{k}^{n}.
\end{equation*}

iii) We use that $D_{k}^{n}D_{k}^{m}=D_{k}^{n+m}$ and we equalize the
elements from the positions $\left( 1,1\right) $. $\Box \medskip $

\textbf{Remark 1.3}. 1) In [Me; 99], relation $i)$ was proved for $k=3$.
Relation \textit{ii)} is a generalization of the Cassini's identity. Indeed,
if we take $k=2$ and $a_{1}=a_{2}=1,$ we obtain that 
\begin{equation*}
d_{n+1}^{2}-d_{n}d_{n+2}=\left\vert 
\begin{array}{cc}
d_{n+1} & d_{n+2} \\ 
d_{n} & d_{n+1}%
\end{array}%
\right\vert =\left\vert 
\begin{array}{cc}
1 & 1 \\ 
1 & 0%
\end{array}%
\right\vert ^{n}=\left( -1\right) ^{n}.
\end{equation*}%
2) In relation iii) from the above proposition, if we take $k=2,$ we obtain

\begin{equation*}
d_{n\text{+}1}d_{m\text{+}1}+a_{2}d_{m}d_{n}=d_{m+n+1}.
\end{equation*}%
If $a_{1}=1$, we get the well know relation for the Fibonacci numbers%
\begin{equation*}
d_{n\text{+}1}d_{m\text{+}1}+d_{m}d_{n}=d_{m+n+1}.
\end{equation*}%
3) In relation iii) from the above proposition, if we take $k=3,$ we obtain 
\begin{equation*}
d_{m\text{+}2}d_{n\text{+}%
2}+a_{2}d_{m+1}d_{n+1}+a_{3}d_{m+1}d_{n}+a_{3}d_{m}d_{n+1}=d_{m\text{+}n%
\text{+}2}\text{.}
\end{equation*}%
\begin{equation*}
\end{equation*}

\bigskip \textbf{2.} \textbf{Pisano period}%
\begin{equation*}
\end{equation*}

For Fibonacci sequence are interesting its properties when it is reduced
modulo $m$. For example, it is well known that this sequence is periodic and
this period is called Pisano's period, denoted by $\pi \left( m\right) $.
For the Fibonacci sequence, $\pi \left( m\right) $ is even, for $m>2$ (see
[Re; 13]). In the following, we generalize this notion for the sequence
generated by the $k-$terms recurrence given in the relation $\left(
1.1\right) $.

First of all, we will prove that the sequence $\left( d_{n}\right) $, given
in $\left( 1.1\right) $, is periodic.\medskip

\textbf{Proposition 2.1}. \textit{The sequence generated by} $\left(
d_{n}\right) $, \textit{considered mod} $m$\textit{, is periodic.\medskip }

\textbf{Prof}. The terms of the sequence $\left( d_{n}\right) $ talked mod $%
m $ can have only the values $\{0,1,2,...,m-1\}$. Since $d_{n}$=$%
a_{1}d_{n-1} $+$a_{2}d_{n-2}$+$...$+$a_{k}d_{n-k}$, we remark that if the
sequence $d_{n-1}$, $d_{n-2}$,$...d_{n-k}$ is repeated mod $m$ from a step,
then the entire sequence is repeated mod $m$. Therefore, there are at most $%
m^{k}$ possible choices for the sequence $d_{n-1}$, $d_{n-2}$,$...d_{n-k}$
and the sequence is periodic.$\Box \medskip $

\textbf{Remark 2.2.} The above proposition generalized Theorem 1 from [Wa;
60] to sequence $\left( d_{n}\right) $.\medskip

\textbf{Definition 2.3.} The \textit{Pisano period} for the sequence $\left(
d_{n}\right) $ is the period with which the sequence $\left( d_{n}\right) $
taken modulo $m$ repeats. We denote this number with $\pi \left( m\right) $%
.\medskip

\textbf{Remark 2.4.} If we consider \ the matrix $D_{k},$ associated to the
difference equation $\left( 1.1\right) ,$ we have that $D_{k}^{\pi \left(
m\right) }\equiv I_{k}$ mod $m,$ where $I_{k}$ is the unity matrix of order $%
k.$ From here, we get 
\begin{equation*}
(\left( -1\right) ^{(k+1)}a_{k})^{\pi \left( m\right) }\equiv 1~\text{mod }m.
\end{equation*}%
\qquad

It results that $ord\left( \left( -1\right) ^{(k+1)}a_{k}\right) \mid \pi
\left( m\right) .\medskip $

\textbf{Theorem 2.5.} \textit{Withe the above notations, the following
statements are true:}

\textit{i)} \textit{If} $s_{1}\mid s_{2}$\textit{, then} $\pi \left(
s_{1}\right) \mid \pi \left( s_{2}\right) $.

\textit{ii) }$\pi \left( \left[ s_{1},s_{2}\right] \right) =\left[ \pi
\left( s_{1}\right) ,\pi \left( s_{2}\right) \right] ,$ \textit{where} $%
s_{1},s_{2}$ \textit{are positive integers and}\newline
$\left[ s_{1},s_{2}\right] =~lcd\left( s_{1},s_{2}\right) $\textit{, the
least common multiple.}

\textit{ii)} $\pi \left( p^{r+1}\right) =$ $\pi \left( p^{r}\right) $ 
\textit{or} $p\pi \left( p^{r}\right) ,$ \textit{with} $\ p$ \textit{a prime
number and }$r$\textit{\ an integer,} $r\geq 1$.

\textit{iv) If} \ $\pi \left( p^{r}\right) \neq \pi \left( p^{r+1}\right) $%
\textit{, then} $\pi \left( p^{r+1}\right) \neq \pi \left( p^{r+2}\right) ,$ 
\textit{with} $\ p$ \textit{a prime number and }$r$\textit{\ an integer,} $%
r\geq 1$.

\textit{v) If} $D_{k}$ \textit{is a diagonalizable matrix and} $p$ \textit{%
is an odd prime, with} $p\nmid a_{k}$\textit{, then} $\pi \left( p\right)
\mid p-1$.

\textit{vi}) \textit{If \ for the numbers sequence} $d_{n}$=$a_{1}d_{n-1}$+$%
a_{2}d_{n-2}$+$...$+$a_{k}d_{n-k},$\newline
$d_{0}$=$d_{1}$=...=$d_{k-2}$=$0,d_{k-1}$=$1,$ \textit{all~}$a_{i}$ \textit{%
are odd, for } $i\in \{1,2,...,k\}$, \textit{then} $\pi \left( 2\right) =$ $%
k+1.\medskip $

\textbf{Proof.} i) Since $s_{1}\mid s_{2},$ we have $D_{k}^{\pi \left(
s_{2}\right) }\equiv I_{k}$ \textit{mod} $s_{1},$ therefore $\pi \left(
s_{1}\right) \mid \pi \left( s_{2}\right) .$

ii) Let $s=\left[ s_{1},s_{2}\right] $. We have $D_{k}^{\pi \left( s\right)
}\equiv I_{k}$ \textit{mod} $s,$ therefore $D_{k}^{\pi \left( s\right)
}\equiv I_{k}$ \textit{mod} $s_{1}$ and $D_{k}^{\pi \left( s\right) }\equiv
I_{k}$ \textit{mod} $s_{2}.$ It results $\pi \left( s_{1}\right) \mid \pi
\left( s\right) $, $\pi \left( s_{2}\right) \mid \pi \left( s\right) $ and
from here,we have $\left[ \pi \left( s_{1}\right) ,\pi \left( s_{2}\right) %
\right] \mid \pi \left( s\right) .~$For the converse, \ we know that $%
D_{k}^{\pi \left( s_{1}\right) }\equiv I_{k}$ \textit{mod} $s_{1}$ and $%
D_{k}^{\pi \left( s_{2}\right) }\equiv I_{k}$ \textit{mod} $s_{2}$,
therefore $D_{k}^{\left[ \pi \left( s_{1}\right) ,\pi \left( s_{2}\right) %
\right] }\equiv I_{k}$ \textit{mod} $s_{1}$ and $D_{k}^{\left[ \pi \left(
s_{1}\right) ,\pi \left( s_{2}\right) \right] }\equiv I_{k}$ \textit{mod} $%
s_{2}$. From here, we obtain $D_{k}^{\left[ \pi \left( s_{1}\right) ,\pi
\left( s_{2}\right) \right] }\equiv I_{k}$ \textit{mod} $s$, thus $\pi
\left( s\right) \mid \left[ \pi \left( s_{1}\right) ,\pi \left( s_{2}\right) %
\right] $.

iii) Supposing that $D_{k}^{\pi \left( p^{r}\right) }\equiv I_{k}$ \textit{%
mod} $p^{r},$ we have $D_{k}^{\pi \left( p^{r}\right) }=I_{k}+p^{r}A$, with $%
A$ a matrix of order $k$. Therefore, $D_{k}^{p\pi \left( p^{r}\right)
}=(I_{k}+p^{r}A)^{p}$ in which, except first term $I_{k}$, all term are
divisible with $p^{r+1}$. It results that $D_{k}^{p\pi \left( p^{r}\right)
}\equiv I_{k}$ \textit{mod} $p^{r+1}$. From here, we obtain that $\pi \left(
p^{r+1}\right) \mid p\pi \left( p^{r}\right) $, thus, from 1), we get $\pi
\left( p^{r+1}\right) =$ $\pi \left( p^{r}\right) $ or $\pi \left(
p^{r+1}\right) =$ $p\pi \left( p^{r}\right) .$

iv) If \ $\pi \left( p^{r}\right) \neq \pi \left( p^{r+1}\right) $\textit{, }%
then, from 3), we have that $\pi \left( p^{r+1}\right) =p\pi \left(
p^{r}\right) $ and $\pi \left( p^{r+2}\right) =$ $\pi \left( p^{r+1}\right) $
or $p\pi \left( p^{r+1}\right) $. Since $D_{k}^{p\pi \left( p^{r}\right)
}=(I_{k}+p^{r}A)^{p}$, it results that $D_{k}^{p\pi \left( p^{r}\right) }$
is not congruent with $I_{k}$ \textit{mod} $p^{r+2},$ therefore $\pi \left(
p^{r+1}\right) \neq \pi \left( p^{r+2}\right) $.

v) \ Supposing that $D_{k}$ is diagonalizable, there are the matrix $A$ and $%
\,B$ of order $k$ such that $B=A^{-1}D_{k}A,$ with $A$ an invertible matrix
and 
\begin{equation*}
B=\left( 
\begin{array}{ccccc}
\lambda _{1} & 0 & 0 & ... & 0 \\ 
0 & \lambda _{2} & 0 & ... & 0 \\ 
0 & 0 & \lambda _{3} & ... & 0 \\ 
... & ... & ... & ... & ... \\ 
0 & 0 & ... & 0 & \lambda _{k}%
\end{array}%
\right) .
\end{equation*}

Applying the Fermat's little Theorem, we have $B^{p-1}\equiv I_{k}$ \textit{%
mod} $p$, therefore $D_{k}^{p-1}\equiv I_{k}$ \textit{mod} $p$ thus $\pi
\left( p\right) \mid p-1$.

vi) Let $d_{n}$=$a_{1}d_{n-1}$+$a_{2}d_{n-2}$+$...$+$a_{k}d_{n-k},$\newline
$d_{0}$=$d_{1}$=...=$d_{k-2}$=$0,d_{k-1}$=$1$. If we consider this sequence
modulo $2$, we obtain $d_{n}$=$d_{n-1}$+$d_{n-2}$+$...$+$d_{n-k},d_{0}$=$%
d_{1}$=...=$d_{k-2}$=$0,d_{k-1}$=$1$. Therefore, we get the following
sequence modulo $2$:

-first $k+1$ terms $d_{0}$, $d_{1}$,..., $d_{k-2}$, $d_{k-1},d_{k}$ are: $%
0,0,...,0,1,1;$

-the next $k+1$ terms $d_{k+1}$, $d_{k+2}$,..., $d_{2k-1}$, $d_{2k},d_{2k+1}$
are: $0,0,...,0,1,1$. From here, it is clear that $\pi \left( 2\right) =k+1$.%
$\Box \medskip $

\textbf{Remark 2.6. }The above theorem generalized for difference equation $%
\left( 1.1\right) \,$\ some results obtained in [Re; 13], Theorem 1,
Proposition 1, Proposition 2 and Theorem 3.\medskip

\textbf{Example 2.7.}

i) Let $k=3$ and $d_{n}$=$a_{1}d_{n-1}$+$a_{2}d_{n-2}$+$a_{3}d_{n-3},d_{0}$=$%
d_{1}$=$0,d_{2}$=$1.$ Supposing that $a_{1},a_{2},a_{3}$ are odd, we
consider this difference equation modulo $2$. It results $d_{n}$=$d_{n-1}$+$%
d_{n-2}$+$d_{n-3},d_{0}$=$d_{1}$=$0,d_{2}$=$1$ and we obtain the sequence:%
\newline
$d_{0}$, $d_{1}$,$d_{2}$, $%
d_{3},d_{4},d_{5},d_{6},d_{7},d_{8},d_{9},d_{10},d_{11}$, etc$:$\newline
$\mathbf{0,0,1,1},0,0,1,1,\mathbf{0,0,1,1},$etc.\newline
Therefore, $\pi \left( 2\right) =k+1=4.$

ii) For the case when some coefficients $a_{i}$ are even, it is difficult to
compute $\pi \left( 2\right) $ since it depends of the number and position
of even coefficients. We give an example for $k=3$, $a_{1},a_{3}$ odd and $%
a_{2}$ even. We have the equation $d_{n}$=$d_{n-1}$+$d_{n-3},d_{0}$=$d_{1}$=$%
0,d_{2}$=$1$ and we obtain the sequence:\newline
$d_{0}$, $d_{1}$,$d_{2}$, $%
d_{3},d_{4},d_{5},d_{6},d_{7},d_{8},d_{9},d_{10},d_{11},d_{12},d_{13},$etc$:$%
\newline
$\mathbf{0,0,1,1,1,0,1},0,0,1,1,1,0,1,0$, etc. Therefore, $\pi \left(
2\right) =7.$

iii) Let $a_{1}=4,a_{2}=-5,a_{3}=2,p=3$. Therefore, we get the matrix%
\begin{equation*}
D_{3}=\left( 
\begin{array}{ccc}
4 & -5 & 2 \\ 
1 & 0 & 0 \\ 
0 & 1 & 0%
\end{array}%
\right) ,
\end{equation*}%
and the associated difference equation $d_{n}$=$%
4d_{n-1}-5d_{n-2}+2d_{n-3},d_{0}$=$d_{1}$=$0,d_{2}$=$1$. For this difference
equation of degree three, we will compute $\pi \left( 3\right) $ \textit{%
modulo} $3$. It results $d_{0}$=$d_{1}$=$0,d_{2}$=$1,d_{3}=1,d_{4}=2,$%
\newline
$d_{5}=2,d_{6}=0,d_{7}=0,d_{8}=1,etc$. Therefore $\pi \left( 3\right) =6.$

Now, we work in $\mathbb{Z}_{9}$ and we will compute $\pi \left(
3^{2}\right) =\pi \left( 9\right) .$ We obtain $d_{0}$=$d_{1}$=$0,d_{2}$=$%
1,d_{3}=1,d_{4}=2,d_{5}=5,d_{6}=0,d_{7}=0,d_{8}=1,etc$. Therefore $\pi
\left( 9\right) =6.$

From Theorem 2.5, iii), it results that $\pi \left( 27\right) =\pi \left(
3\right) $ or $\pi \left( 27\right) =3\cdot \pi \left( 3\right) =18.$
Working in $\mathbb{Z}_{27},$ we get $d_{0}$=$d_{1}$=$0,d_{2}$=$%
1,d_{3}=1,d_{4}=2,d_{5}=5,d_{6}=9,...,$ therefore $\pi \left( 27\right) =18$.

iv) We work on $\mathbb{Z}_{5}$, therefore $p=5$. Let $%
a_{1}=6,a_{2}=-11,a_{3}=6$. Therefore, we get the matrix%
\begin{equation*}
D_{3}=\left( 
\begin{array}{ccc}
6 & -11 & 6 \\ 
1 & 0 & 0 \\ 
0 & 1 & 0%
\end{array}%
\right)
\end{equation*}
and the associated difference equation $d_{n}$=$%
6d_{n-1}-11d_{n-2}+6d_{n-3},d_{0}$=$d_{1}$=$0,d_{2}$=$1$. For this
difference equation of degree three, we will compute $\pi \left( 5\right) $ 
\textit{modulo} $5$. We have that the eigenvalues of the matrix $D_{3}$ are $%
\{1,2,3\}$, therefore $D_{3}$ has the diagonal form 
\begin{equation*}
Diag=\left( 
\begin{array}{ccc}
1 & 0 & 0 \\ 
0 & 2 & 0 \\ 
0 & 0 & 3%
\end{array}%
\right) \text{.}
\end{equation*}
It results that $\pi \left( 5\right) \mid 4$. We compute the terms of the
sequence $\left( d_{n}\right) _{n\in \mathbb{N}}$ and we obtain $d_{0}$=$%
d_{1}$=$0,d_{2}$=$1,d_{3}=1,d_{4}=0,d_{5}=0,d_{6}=1$, therefore, in this
situation, $\pi \left( 5\right) =4$. 
\begin{equation*}
\end{equation*}

\bigskip \textbf{3. An application in Cryptography}%
\begin{equation*}
\end{equation*}

In the following, we will give an application of difference equation in
Cryptography. In [Ko; 94] was presented enciphering matrices. Since these
matrices must be invertible, we will use as enciphering matrices the
matrices of the form $D_{k}^{n}$, since the determinant of such a matrix is
known, namely \textit{det}$D_{k}=\left( -1\right) ^{k+1}a_{k}$ (see
[Jo]).\medskip

\textbf{Encryption\medskip }

1) Let $\mathcal{A}$ be an alphabet with $N$ letters, labelled from $0$ to $%
N-1$, let $m$ be a plain text and let $q$ be a number obtained by using the
labels of the letters from the plain text $m$. We chose $k$ a nonzero
natural number, we split the text $m$ in blocks of length $k$, $%
m_{1},m_{2},...,m_{r}$. To these blocks correspond the vectors $%
v_{1},v_{2},...,v_{r}$ of dimension $k$.

2) We chose arbitrary integers,\ $a_{1},...,a_{k}\in
\{0,1,2,...,N-1\},a_{k}\neq 0~$and the general $k-$terms recurrence equation
given by the relation $\left( 1.1\right) $, $n,k\in \mathbb{N},k\geq 2,n\geq
k$.

3) We chose the matrix $D_{k}$ associated to the sequence $a_{1},...,a_{k}$
and given by the relation $\left( 1.2\right) $.

4) We chose $n\in \mathbb{N}-\{0\}$ the power of the matrix $D_{k}.$

5) We use as an enciphering matrice the matrix $D_{k}^{n}$, therefore the
encryption key will be $\left( k,N,a_{1},...,a_{k},n\right) .$

6) We work on $\mathbb{Z}_{\mathbb{N}},$therefore, the encrypted text is
given by the following formula%
\begin{equation}
C=D_{k}^{n}V,  \tag{3.1.}
\end{equation}%
where $V=\left( v_{1}v_{2}...v_{r}\right) \in \mathcal{M}_{k\times r}\left( 
\mathbb{Z}_{N}\right) $ is a matrix with columns $v_{1},v_{2},...,v_{r}$ and 
$C\in \mathcal{M}_{k\times r}\left( \mathbb{Z}_{N}\right) $ is a matrix with
columns $\left( c_{1},c_{2},...,c_{r}\right) $.

We remark that it is very important to know the Pisano period $\pi \left(
N\right) .$ In this case, the key can have the form 
\begin{equation*}
\left( k,N,a_{1},...,a_{k},n\right) =\left( k,N,a_{1},...,a_{k},l\pi \left(
N\right) +n\right) ,l\in \mathbb{Z}
\end{equation*}%
and the encryption and decryption algorithms are the same since 
\begin{equation*}
D_{k}^{n}=D_{k}^{l\pi \left( N\right) +n},l\in \mathbb{Z}\text{.}
\end{equation*}%
In this way, the key can be made a litle bit hard. Therefore, to make the
algorithm faster, we can choose $n\in \{1,2,...,\pi \left( N\right) -1\}$%
.\medskip

\textbf{Decryption\medskip }

-The decryption key is the same as the encryption key and will be $\left(
k,N,a_{1},...,a_{k},n\right) .$

-We obtain the matrices $D_{k},D_{k}^{n}$ and $D_{k}^{-n}.$

-The decrypted text is 
\begin{equation}
V=D_{k}^{-n}C,  \tag{3.2.}
\end{equation}%
with \ $V=\left( v_{1}v_{2}...v_{r}\right) \in \mathcal{M}_{k\times r}\left( 
\mathbb{Z}_{N}\right) $.\medskip\ 

\textbf{Example 3.1.\medskip \medskip }

We work in $\mathbb{Z}_{2}$ and we consider the encryption key $\left(
3,2,1,1,1,3\right) $. Therefore, we get the matrix%
\begin{equation*}
D_{3}=\left( 
\begin{array}{ccc}
1 & 1 & 1 \\ 
1 & 0 & 0 \\ 
0 & 1 & 0%
\end{array}%
\right) .
\end{equation*}%
Since 
\begin{equation*}
D_{3}^{n}=\left( 
\begin{array}{ccc}
d_{n+2} & bd_{n+1}+cd_{n} & cd_{n+1} \\ 
d_{n+1} & bd_{n}+cd_{n-1} & cd_{n} \\ 
d_{n} & bd_{n-1}+cd_{n-2} & cd_{n-1}%
\end{array}%
\right) .
\end{equation*}%
Using Example 2.7, i), we get the matrix%
\begin{equation*}
D_{3}^{3}=\left( 
\begin{array}{ccc}
d_{5} & d_{4}+d_{3} & d_{4} \\ 
d_{4} & d_{3}+d_{2} & d_{3} \\ 
d_{3} & d_{2}+d_{1} & d_{2}%
\end{array}%
\right) =\left( 
\begin{array}{ccc}
0 & 1 & 0 \\ 
0 & 0 & 1 \\ 
1 & 1 & 1%
\end{array}%
\right) .
\end{equation*}

We consider the alphabet $\mathcal{A}=\{A,B\}$, labeled $\{0,1\}$ and the
plain text $m=ABBAAB$. The number obtained by using the labels of the
letters from the plain text $m$ is $q=011001$. The obtained vectors are $%
v_{1}=\left( 
\begin{array}{c}
0 \\ 
1 \\ 
1%
\end{array}%
\right) ,v_{2}=\left( 
\begin{array}{c}
0 \\ 
0 \\ 
1%
\end{array}%
\right) .$ Therefore, the encrypted text is 
\begin{equation*}
C=\left( 
\begin{array}{ccc}
0 & 1 & 0 \\ 
0 & 0 & 1 \\ 
1 & 1 & 1%
\end{array}%
\right) \left( 
\begin{array}{cc}
0 & 0 \\ 
1 & 0 \\ 
1 & 1%
\end{array}%
\right) =\left( 
\begin{array}{cc}
1 & 0 \\ 
1 & 1 \\ 
0 & 1%
\end{array}%
\right) .
\end{equation*}%
Thus, \thinspace $c_{1}=\left( 
\begin{array}{c}
1 \\ 
1 \\ 
0%
\end{array}%
\right) ,c_{2}=\left( 
\begin{array}{c}
0 \\ 
1 \\ 
1%
\end{array}%
\right) $ and the obtained encrypted text is $BBAABB.$

We remark that the value $n$ can be bigger than $3$, for example $31$, which
can help us to have a little bit hard key. But, since $31$\textit{\ mod} $%
\pi \left( 2\right) =31$\textit{mod} $4=3$, this can help us in the
encrypting and decrypting process. Therefore, the value $n$ can be taken
from the set$~\{1,2,...,\pi \left( N\right) -1\}.$

To decrypt the text, we use the same key $\left( 3,2,1,1,1,3\right) $ and we
obtain 
\begin{equation*}
V=D_{3}^{-3}C=\left( 
\begin{array}{ccc}
0 & 1 & 0 \\ 
0 & 0 & 1 \\ 
1 & 1 & 1%
\end{array}%
\right) ^{-1}\left( 
\begin{array}{cc}
1 & 0 \\ 
1 & 1 \\ 
0 & 1%
\end{array}%
\right) =
\end{equation*}%
\begin{equation*}
=\left( 
\begin{array}{ccc}
1 & 1 & 1 \\ 
1 & 0 & 0 \\ 
0 & 1 & 0%
\end{array}%
\right) \left( 
\begin{array}{cc}
1 & 0 \\ 
1 & 1 \\ 
0 & 1%
\end{array}%
\right) =\left( 
\begin{array}{cc}
0 & 0 \\ 
1 & 0 \\ 
1 & 1%
\end{array}%
\right) ,
\end{equation*}%
where $V=\left( v_{1}v_{2}...v_{r}\right) \in \mathcal{M}_{k\times r}\left( 
\mathbb{Z}_{\mathbb{N}}\right) $ is a matrix with columns $%
v_{1},v_{2},...,v_{r}$, that means we get the initial message $ABBAAB$%
.\medskip

\textbf{Example 3.2.\medskip\ }

We use Example 2.7, iii). For $a_{1}=4,a_{2}=-5,a_{3}=2,p=3$,we consider the
matrix%
\begin{equation*}
D_{3}=\left( 
\begin{array}{ccc}
4 & -5 & 2 \\ 
1 & 0 & 0 \\ 
0 & 1 & 0%
\end{array}%
\right)
\end{equation*}%
and the associated difference equation $d_{n}$=$%
4d_{n-1}-5d_{n-2}+2d_{n-3},d_{0}$=$d_{1}$=$0,d_{2}$=$1$.

We consider the alphabet $\mathcal{A}=\{A,B,C,...,Z,\ast \}$, labeled $%
\{0,1,2,...,26\},$ where "$\ast $" is the blank character and the plain text 
$m=SUCCESS\ast \ast $. Therefore, we work in $\mathbb{Z}_{27}$ and, from
Example 2.7, iii), we know that $\pi \left( 27\right) =18$. We consider the
encryption key $\left( 3,27,4,-5,2,2\right) $. The number obtained by using
the labels of the letters from the plain text $m$ is $q=182002020418182626$.
The obtained vectors are $v_{1}=\left( 
\begin{array}{c}
18 \\ 
20 \\ 
02%
\end{array}%
\right) ,v_{2}=\left( 
\begin{array}{c}
02 \\ 
04 \\ 
18%
\end{array}%
\right) ,v_{2}=\left( 
\begin{array}{c}
18 \\ 
26 \\ 
26%
\end{array}%
\right) $. We have

\begin{equation*}
D_{3}^{2}=\left( 
\begin{array}{ccc}
2 & 24 & 2 \\ 
1 & 22 & 2 \\ 
1 & 0 & 0%
\end{array}%
\right) .
\end{equation*}%
Therefore, the encrypted text is 
\begin{equation*}
C\text{=}\left( 
\begin{array}{ccc}
2 & 24 & 2 \\ 
1 & 22 & 2 \\ 
1 & 0 & 0%
\end{array}%
\right) \left( 
\begin{array}{ccc}
18 & 2 & 18 \\ 
20 & 4 & 26 \\ 
2 & 18 & 26%
\end{array}%
\right) \text{=}
\end{equation*}
\begin{equation*}
\text{=}\left( 
\begin{array}{ccc}
520 & 136 & 712 \\ 
462 & 126 & 642 \\ 
18 & 2 & 18%
\end{array}%
\right) \text{ \textit{mod }}27\text{=}\left( 
\begin{array}{ccc}
7 & 1 & 10 \\ 
3 & 18 & 21 \\ 
18 & 2 & 18%
\end{array}%
\right) .
\end{equation*}

Therefore, the encrypted text is $HDSBSCKVS$. We remark that the analysis of
the frequency of the letters in the text can't be applied here, since the
same letters are encrypted in different characters and different characters
can be encrypted in the same letter, as we can see in this example. In this
way, this new method provides multiple ways of finding the encryption and
decryption keys, thus the obtained encrypted texts are hard to break.

The decryption key is the same as the encryption key, namely $\left(
3,27,4,-5,2,2\right) $.

We obtain 
\begin{equation*}
V=D_{3}^{-2}C=\left( 
\begin{array}{ccc}
2 & 24 & 2 \\ 
1 & 22 & 2 \\ 
1 & 0 & 0%
\end{array}%
\right) ^{-1}\left( 
\begin{array}{ccc}
7 & 1 & 10 \\ 
3 & 18 & 21 \\ 
18 & 2 & 18%
\end{array}%
\right) =
\end{equation*}%
\begin{equation*}
\left( 
\begin{array}{ccc}
0 & 0 & 1 \\ 
14 & 13 & 13 \\ 
8 & 6 & 5%
\end{array}%
\right) \left( 
\begin{array}{ccc}
7 & 1 & 10 \\ 
3 & 18 & 21 \\ 
18 & 2 & 18%
\end{array}%
\right) =
\end{equation*}%
\begin{equation*}
=\left( 
\begin{array}{ccc}
18 & 2 & 18 \\ 
371 & 274 & 647 \\ 
164 & 126 & 296%
\end{array}%
\right) =\left( 
\begin{array}{ccc}
18 & 2 & 18 \\ 
20 & 4 & 26 \\ 
2 & 18 & 26%
\end{array}%
\right) .
\end{equation*}

We remark that the same result can be obtained in the case when the key is
under the form $\left( 3,27,4,-5,2,l\pi \left( 27\right) +2\right) $, where $%
l\in \mathbb{Z}$.\medskip

\textbf{Remark 3.3.} \ The cryptosystem described in this section is a
symmetric cryptosystem. A problem of such a system is the key transmission.

A solution is to use hybrid cryptosystems (using in the same time symmetric
and asymmetric cryptosystems). In this situation, data can be protected
using a symmetric key, and the symmetric key can be encrypted and
distributed by using a public key.

If are used only a symmetric cryptosystem, the key must be encrypted by
using a key-encrypting keys (KEKs). In this situation, the key-encrypting
keys must be distributed by using an arboreal algorithm. Usually, such an
algorithm must have some steps. These keys (KEKs) are used to encrypt other
keys. The most difficult part is to initialize this process. For initialize
this process, the first key is interchanged between users and this key is
used to encrypt another key, called its successor and so on. The security of
the first exchange step must be high, otherwise all successors keys are
compromised. For the initialization phase of the process, can be used some
procedures of the keys distribution as for example: transmission of the keys
by using certain channels or "face-to-face identification", (see [EPC; 18]).

\begin{equation*}
\end{equation*}

\textbf{4.} \textbf{Applications of some special number sequences and
quaternion elements}

\begin{equation*}
\end{equation*}

Let $\mathbb{H}\left( \alpha ,\beta \right) $ be the generalized\ quaternion
algebra over an arbitrary field $\mathbb{K}$, i.e. the algebra of the
elements of the form $a=a_{1}\cdot 1+a_{2}e_{2}+a_{3}e_{3}+a_{4}e_{4},$
where $a_{i}\in \mathbb{K},i\in \{1,2,3,4\}$, and the elements of the basis $%
\{1,e_{2},e_{3},e_{4}\}$ satisfying the following rules, given in the below
multiplication table:\medskip\ \vspace{3mm}

\begin{center}
\begin{tabular}{c|cccc}
$\cdot $ & $1$ & $e_{2}$ & $e_{3}$ & $e_{4}$ \\ \hline
$1$ & $1$ & $e_{2}$ & $e_{3}$ & $e_{4}$ \\ 
$e_{2}$ & $e_{2}$ & $\alpha $ & $e_{4}$ & $\alpha e_{3}$ \\ 
$e_{3}$ & $e_{3}$ & $-e_{4}$ & $\beta $ & $-\beta e_{2}$ \\ 
$e_{4}$ & $e_{4}$ & $-\alpha e_{3}$ & $\beta e_{2}$ & $-\alpha \beta $%
\end{tabular}%
.\medskip
\end{center}

Let $\overline{a}=a_{1}\cdot 1-a_{2}e_{2}-a_{3}e_{3}-a_{4}e_{4}$ be the
conjugate of the quaternion $a.$ The norm of $a$ is $\boldsymbol{n}\left(
a\right) =a\cdot \overline{a}=a_{1}^{2}-\alpha a_{2}^{2}-\beta
a_{3}^{2}+\alpha \beta a_{4}^{2}$ and the trace of the element $\ a$ is $%
\mathbf{t}\left( a\right) =a+\overline{a}.$\newline
If, \thinspace for $x\in \mathbb{H}\left( \alpha ,\beta \right) $, the
relation $\mathbf{n}\left( x\right) =0$ implies $x=0$, then the algebra $%
\mathbb{H}\left( \alpha ,\beta \right) $ is called a \textit{division}
algebra, otherwise the quaternion algebra is called a \textit{split} algebra.

If $p$ is a prime number, it is known that the quaternion algebra $\mathbb{H}%
_{\mathbb{Z}_{p}}\left( -1,-1\right) $ splits. In the paper [Sa; 17], the
second author determined the Fibonacci quaternions which are zero divisors
in the quaternion algebra $\mathbb{H}_{\mathbb{Z}_{p}}\left( -1,-1\right) ,$
respectively the Fibonacci quaternions which are invertible in the
quaternion algebra $\mathbb{H}_{\mathbb{Z}_{p}}\left( -1,-1\right) .$ 
\newline

In the paper [Gr, Mi, Ma; 15], Grau, Miguel and Oller-Marcen have studied
the quaternion algebra $\mathbb{H}\left( -1,-1\right) $ over a finite
commutative unitary ring $\left( \mathbb{Z}_{n},+,\cdot \right) ,$ where $n$
is a positive integer, $n\geq 3.$\newline

In the following, when $l$ is an odd prime number, we study $l$-quaternions
in quaternion algebra $\mathbb{H}_{\mathbb{Z}_{l}}\left( -1,-1\right) $ and
also we study $l$-quaternions in the quaternion ring $\mathbb{H}_{\mathbb{Z}%
_{l^{r}}}\left( -1,-1\right) ,$ where $r$ is a positive integer, $r\geq 2$.

Let $l$ be a nonzero positive integer. In [Sa; 19], the second author
considered the sequence $\left( a_{n}\right) _{n\geq 0},$ 
\begin{equation*}
a_{n}=la_{n-1}+a_{n-2},\;n\geq 2,a_{0}=0,a_{1}=1
\end{equation*}%
We call these numbers $\left( l,1,0,1\right) $ numbers or $l-$ numbers. We
remark that for $l=1,$ it is obtained the Fibonacci numbers and for $l=2,$
it is obtained the Pell numbers. In the following, we present some
properties of these numbers.\medskip \smallskip \newline

\textbf{Remark 4.1.} ([Sa; 19]). \ Let $\left( a_{n}\right) _{n\geq 0}$ be
the sequence previously defined. Then, the following relations are true:%
\newline
i)%
\begin{equation*}
a_{n}^{2}+a_{n+1}^{2}=a_{2n+1},~\text{for all }n\in \mathbb{N}.
\end{equation*}%
ii) For $\alpha =\frac{l+\sqrt{l^{2}+4}}{2}$ and $\beta =\frac{l-\sqrt{%
l^{2}+4}}{2},$ we obtain that \smallskip\ 
\begin{equation*}
a_{n}=\frac{\alpha ^{n}-\beta ^{n}}{\alpha -\beta }=\frac{\alpha ^{n}-\beta
^{n}}{\sqrt{l^{2}+4}},\ \text{for all }n\in \mathbb{N},
\end{equation*}%
called the \textbf{Binet's formula for the sequence} $\left( a_{n}\right)
_{n\geq 0}.\medskip $\newline

\textbf{Proposition 4.2.} ([Fl, Sa; 19]) \textit{\ Let }$\left( a_{n}\right)
_{n\geq 0}$\textit{\ be the sequence previously defined. The following
relations hold:}

i) \textit{If} $d\mid n,$ \textit{then } $a_{d}\mid a_{n}.$

ii) $a_{m+n}=a_{m}a_{n+1}+a_{m-1}a_{n}.\Box $\smallskip \medskip 

\textbf{Proposition 4.3.} ([Fl, Sa; 19]) \textit{Let} $(a_{n})_{n\geq 0}$ 
\textit{be the sequence previously defined.} \textit{Then, the following
relations are true:}\newline
i) 
\begin{equation*}
a_{n}+a_{n+4}=\left( l^{2}+2\right) a_{n+2};
\end{equation*}%
ii) 
\begin{equation*}
a_{n}+a_{n+8}=\left[ \left( l^{2}+2\right) ^{2}-2\right] a_{n+4};
\end{equation*}%
iii) 
\begin{equation*}
a_{n}+a_{n+2^{k}}=M_{k}a_{n+2^{k-1}},k\geq 3,
\end{equation*}%
\textit{where} $M_{k}=\underset{k-3\;times\;of\;-2}{\underbrace{\left[
\left( \left( \left( l^{2}+2\right) ^{2}-2\right) ^{2}...-2\right) ^{2}-2%
\right] }}.$ \textit{The sequence} $\left( M_{k}\right) _{k\geq 2}$ \textit{%
satisfies the recurrence} $M_{k+1}=M_{k}^{2}-2,$ \textit{for all} $k$$\in $$%
\mathbb{N},$ $k\geq 2,$ $M_{2}=l^{2}+2.\Box $\newline
\smallskip \newline
Using these results, in the following, we give new properties of $l-$
numbers.\newline
We begin with some examples: $a_{0}=0\equiv0$ (mod $l^{2}) ,$\textbf{\ }$%
a_{1}=1$$\equiv $$1$ (mod $l^{2}$), $a_{2}=l\equiv l$ (mod $l^{2}),$ $%
a_{3}=l^{2}+1$$\equiv $$1$ (mod $l^{2}$), $a_{4}=l\left( l^{2}+2\right)
\equiv 2l$ (mod $l^{2}),$ $a_{5}=l^{2}\left( l^{2}+3\right) +1$$\equiv $$1$
(mod $l^{2})$.\medskip\ \smallskip \newline

\textbf{Proposition 4.4.} \textit{Let} $l$ \textit{be a nonzero positive
integer and let} $\left( a_{n}\right) _{n\geq 0}$ \textit{be the sequence of 
} $l-$ \textit{numbers. Then, the following relations hold}:\newline
i) $l\mid a_{n}$ \textit{if and only if} $n$ \textit{is an even number};%
\newline
ii) $a_{n}$$\equiv $$1$ (\textit{mod} $l^{2}$) \textit{if and only if} $n$ 
\textit{is an odd number}.\medskip\ \smallskip 

\textbf{Proof.} The proof is a straightforward calculation, using
Proposition 4.2 (i) and a mathematical induction after $n$$\in $$\mathbb{N}%
.\Box \medskip $\newline

\textbf{Proposition 4.5.} \textit{Let} $l$ \textit{be a nonzero positive
integer and let} $\left( a_{n}\right) _{n\geq 0}$ \textit{be the sequence of 
} $l-$ \textit{numbers. Then, the set} 
\begin{equation*}
M=\left\{ \alpha a_{2n}|n\in \mathbb{N},\alpha \in \mathbb{Z}\right\}
\end{equation*}%
\textit{is a commutative nonunitary ring, with addition and
multiplication.\medskip }

\textbf{Proof.} First remark is that $a_{0}=0$$\in $$M$ is the identity
element for addition on $M$ and $a_{1}=1$$\notin $$M$. It is clear that the
addition and multiplication on $M$ are commutative.\newline
\qquad We prove that $\left( M,+,\cdot \right) $ is a subring of the ring $%
\left( \mathbb{Z},+,\cdot \right) $. According to Proposition 4.4 (i), it
results that there is $\lambda $$\in \mathbb{Z}$ such that $a_{2n}=\lambda
a_{2}$. From here, we obtain that for each $n,m$$\in $$\mathbb{N}$ and for
each $\alpha ,\beta $$\in $$\mathbb{Z},$ we have 
\begin{equation*}
\alpha a_{2n}-\beta a_{2m}=\gamma a_{2r}
\end{equation*}%
and 
\begin{equation*}
\alpha a_{2n}\cdot \beta a_{2m}=\delta a_{2s},
\end{equation*}%
where $r,s$$\in $$\mathbb{N}$ and $\gamma ,\delta $$\in $$\mathbb{Z}.$%
\newline
Therefore, it results that $\left( M,+,\cdot \right) $ is a commutative
nonunitary ring.\newline

The above proposition generalized to a difference equation of degree $k,$
Proposition 3.2 from [FSZ; 19].\medskip

\textbf{Proposition 4.6.} \textit{Let} $l$ \textit{be a nonzero positive
integer and let} $\left( a_{n}\right) _{n\geq 0}$ \textit{be the sequence of}
$l-$ \textit{numbers. Then, the set} 
\begin{equation*}
M=\left\{ \alpha a_{2n}|n\in \mathbb{N},\alpha \in \mathbb{Z}\right\}
\end{equation*}%
\textit{is an ideal of the ring} $\left( \mathbb{Z},+,\cdot \right) ,$ 
\textit{with the property} $M=$ $l\mathbb{Z}$.\medskip \smallskip

\textbf{Proof.} Applying Proposition 4.5 and the fact that $\beta \left(
\alpha a_{2n}\right) \gamma =\alpha \beta \gamma a_{2n}$$\in $$M$, for each $%
n$$\in $$\mathbb{N}$ and for each $\alpha ,\beta ,\gamma $$\in $$\mathbb{Z},$
it results that $M$ is a bilateral ideal of the ring $\left( \mathbb{Z}%
,+,\cdot \right) .$ According to Proposition 4.4 (i) it results that $M=$ $l%
\mathbb{Z}.\Box \medskip $\smallskip \newline
\qquad We consider now the $l-$ numbers $\left( a_{n}\right) _{n\geq 0}$ in
the case when $l$ is an odd prime number. Let $\left( \mathbb{Z}_{l},+,\cdot
\right) $ be a finite field and let $\mathbb{H}_{\mathbb{Z}_{l}}\left(
-1,-1\right) $ be the quaternion algebra$.$ It is known that this quaternion
algebra splits. Let $\left\{ 1,e_{1},e_{2},e_{2},e_{3}\right\} $ be a basis
of this algebra and let $A_{n}$ be the $n$th $l-$ quaternion, 
\begin{equation*}
A_{n}=a_{n}1+a_{n+1}e_{1}+a_{n+2}e_{2}+a_{n+3}e_{3}.
\end{equation*}%
\qquad \qquad 

In the following, we will determine the invertible $l-$ quaternions from the
quaternion algebra $\mathbb{H}_{\mathbb{Z}_{l}}\left( -1,-1\right) $.\medskip

\textbf{Proposition 4.7.} \textit{Let} $l$ \textit{be an odd prime integer,
let} $\left( a_{n}\right) _{n\geq 0}$ \textit{be the sequence of } $l-$ 
\textit{numbers and let } $\mathbb{H}_{\mathbb{Z}_{l}}\left( -1,-1\right) $ 
\textit{be the quaternion algebra}$.$ \textit{Then, all the} $n$\textit{th} $%
l-$ \textit{quaternions are invertible in the quaternion algebra} $\mathbb{H}%
_{\mathbb{Z}_{l}}\left( -1,-1\right) $.\medskip

\textbf{Proof.} In the paper [Fl, Sa; 19], we obtained that the norm of the $%
n$th $l-$ quaternion is $n\left( A_{n}\right) =\left( l^{2}+2\right)
a_{2n+3}.$ Applying Proposition 4.4. (ii), it results that $n\left(
A_{n}\right) $$\equiv $$2$ mod $l$. Since $l$ is odd, we have that $n\left(
A_{n}\right) $$\neq $$\widehat{0}$ in $\mathbb{Z}_{l},$ for all $n$$\in $$%
\mathbb{N}.$ Thus, all the $n$ th $l-$ quaternions are invertible in the
quaternion algebra $\mathbb{H}_{\mathbb{Z}_{l}}\left( -1,-1\right) $.$\Box $%
\smallskip

In the paper [Gr, Mi, Ma; 15] Grau, Miguel and Oller-Marcen proved that when 
$n=p^{r},$ whit $p$ an odd prime and $r$ a positive integer, the quaternion
ring $\mathbb{H}_{\mathbb{Z}_{p^{r}}}\left( -1,-1\right) $ is isomorphic
with the matrix ring $M_{2}\left( \mathbb{Z}_{p^{r}}\right) $ (see [Gr, Mi,
Ma; 15], Proposition 4), therefore the quaternion ring $\mathbb{H}_{\mathbb{Z%
}_{p^{r}}}\left( -1,-1\right) $ splits. We want to determine how many
invertible $n$th $l-$ quaternions are in this quaternion ring, when $p=l$%
.\medskip \smallskip $\Box $\qquad

\textbf{Proposition 4.8.} \textit{Let} $l$ \textit{be an odd prime integer,
let} $\left( a_{n}\right) _{n\geq 0}$ \textit{be the sequence of} $l-$ 
\textit{numbers and let }$\mathbb{H}_{\mathbb{Z}_{l^{r}}}\left( -1,-1\right) 
$ \textit{be the quaternion ring}$.$ \textit{Then, the following statements
are true}:\newline
i) \textit{All the} $n$\textit{th} $l-$ \textit{quaternions are invertible
in the quaternion ring~}$\mathbb{H}_{\mathbb{Z}_{l^{2}}}\left( -1,-1\right) $%
;\newline
ii) \textit{All the} $n$\textit{th} $l-$ \textit{quaternions are invertible
in the quaternion ring} $\mathbb{H}_{\mathbb{Z}_{l^{r}}}\left( -1,-1\right) $%
, \textit{where} $r$ \textit{is a positive integer}, $r\geq 3$.\medskip\ 
\newline
\smallskip \qquad \textbf{Proof.} i) Similar to the proof of Proposition
4.7, applying Proposition 4.4. (ii), it results that $n\left( A_{n}\right) $$%
\equiv $$2$ (mod $l^{2}$). We obtain that all the $n$th $l-$ quaternions are
invertible in the quaternion ring $\mathbb{H}_{\mathbb{Z}_{l^{2}}}\left(
-1,-1\right) $.\newline
ii) Let $r$ be a positive integer, $r\geq 3$. Since $l$ is odd and $n\left(
A_{n}\right) $$\equiv $$2$ (mod $l^{2}$), it results that $n\left(
A_{n}\right) $$\not\equiv $$0$ (mod $l^{r}$). Therefore, all the $n$th $l-$
quaternions are invertible in the quaternion ring $\mathbb{H}_{\mathbb{Z}%
_{l^{r}}}\left( -1,-1\right) $.$\Box $\smallskip 

\textbf{Proposition 4.9.} \textit{Let} $(a_{n})_{n\geq 0}$ \textit{be the
sequence of} $\ l-$ \textit{numbers and let} $(M_{k})_{k\geq 2}$ \textit{be
the sequence from Proposition 4.3}. \textit{Then, the following relation is
true:}\newline
\begin{equation*}
a_{n}+a_{n+3\cdot 2^{k}}=M_{k}\left( M_{k}^{2}-3\right) a_{n+3\cdot
2^{k-1}},k\geq 2.
\end{equation*}%
\qquad

\textbf{Proof.} Let $n,k$$\in $$\mathbb{N},$ $k\geq 2.$ Applying Proposition
4.3 (iii), we have: 
\begin{equation*}
a_{n}+a_{n+3\cdot 2^{k}}=\left( a_{n}+a_{n+2^{k}}\right) +\left(
a_{n+2^{k}}+a_{n+2^{k+1}}\right) +
\end{equation*}%
\begin{equation*}
+\left( a_{n+2^{k+1}}+a_{n+3\cdot 2^{k}}\right) -2\left(
a_{n+2^{k}}+a_{n+2^{k+1}}\right) =
\end{equation*}%
\begin{equation*}
=M_{k}a_{n+2^{k-1}}+M_{k}a_{n+3\cdot 2^{k-1}}+M_{k}a_{n+5\cdot
2^{k-1}}-2M_{k}a_{n+3\cdot 2^{k-1}}=
\end{equation*}%
\begin{equation*}
=M_{k}a_{n+2^{k-1}}+M_{k}a_{n+5\cdot 2^{k-1}}-M_{k}a_{n+3\cdot 2^{k-1}}=
\end{equation*}%
\begin{equation*}
=M_{k}M_{k+1}a_{n+3\cdot 2^{k-1}}-M_{k}a_{n+3\cdot 2^{k-1}}=
\end{equation*}%
\begin{equation*}
=M_{k}\left( M_{k+1}-1\right) a_{n+3\cdot 2^{k-1}}=M_{k}\left(
M_{k}^{2}-3\right) a_{n+3\cdot 2^{k-1}}.
\end{equation*}%
$\Box $\smallskip \smallskip \newline

\textbf{Proposition 4.10.} \textit{Let} $l$ \textit{be an odd prime number,
let} $\left( a_{n}\right) _{n\geq 0}$ \textit{be the sequence of} $l-$ 
\textit{numbers and let }$\mathbb{H}_{\mathbb{Z}_{l}}\left( -1,-1\right) $ 
\textit{be the quaternion algebra}$.$ \textit{Let} $A_{n}$ \textit{be the} $%
n $\textit{th} $l-$ \textit{quaternion in the quaternion algebra} $\mathbb{H}%
_{\mathbb{Z}_{l}}\left( -1,-1\right) $. \textit{Then, we have}:\newline
\begin{equation*}
A_{n}=A_{n+2},\text{ \textit{for all} }n\in \mathbb{N}\text{.}
\end{equation*}

\textbf{Proof.} We use that $a_{n}$$\equiv $$a_{n+2}$ mod $l$ and
Proposition 4.4. $\Box $\smallskip \smallskip \newline

\textbf{Proposition 4.11.} \textit{Let} $l$ \textit{be an odd prime integer,
let} $\left( a_{n}\right) _{n\geq 0}$ \textit{be the sequence of } $l-$ 
\textit{numbers and let }$\mathbb{H}_{\mathbb{Z}_{l}}\left( -1,-1\right) $ 
\textit{be the quaternion algebra}$.$ \textit{Let} $A_{n}$ \textit{be the} $%
n $\textit{th} $l-$ \textit{quaternion in the quaternion algebra} $\mathbb{H}%
_{\mathbb{Z}_{l}}\left( -1,-1\right) $ \textit{and} $m$ \textit{be a fixed
positive integer. Then, the set} 
\begin{equation*}
M^{^{\prime }}=\left\{ \alpha A_{m+2n}|n\in \mathbb{N},\alpha \in \mathbb{Z}%
\right\} \cup \left\{ 0\right\}
\end{equation*}%
\textit{is a} $\mathbb{Z}-$ \textit{module}.\medskip

\textbf{Proof.} It results from Proposition 4.10. $\Box $\smallskip

\textbf{Proposition 4.12.} \textit{Let} $l$ \textit{be an odd prime integer,
let} $\left( a_{n}\right) _{n\geq 0}$ \textit{be the sequence of} $l-$ 
\textit{numbers and let } $\mathbb{H}_{\mathbb{Z}_{l^{2}}}\left(
-1,-1\right) $ \textit{be the quaternion ring}$.$ \textit{Let} $A_{n}$ 
\textit{be the} $n$\textit{th} $l-$ \textit{quaternion in the quaternion ring%
} $\mathbb{H}_{\mathbb{Z}_{l^{2}}}\left( -1,-1\right) .$ \textit{Let} $%
(M_{k})_{k\geq 2}$ \textit{be the sequence from Proposition 4.3}. \textit{%
Therefore, we have}:\newline
i) 
\begin{equation*}
A_{n}+A_{n+2^{k}}=\widehat{2}A_{n+2^{k-1}},k\geq 2,
\end{equation*}%
ii) 
\begin{equation*}
A_{n}+A_{n+3\cdot 2^{k}}=\widehat{2}A_{n+3\cdot 2^{k-1}},k\geq 2.
\end{equation*}%
iii) 
\begin{equation*}
A_{n}+A_{n+1}+...+A_{n+2l^{2}-1}=\widehat{0}\in \mathbb{Z}_{l^{2}}.
\end{equation*}%
\smallskip \newline

\textbf{Proof.} i) Using Proposition 4.3 and working in $\mathbb{Z}_{l^{2}}$
we obtain that 
\begin{equation*}
A_{n}+A_{n+2^{k}}=\widehat{M_{k}}A_{n+2^{k-1}},k\geq 2.
\end{equation*}%
It is easy to remark that $M_{k}\equiv 2$ mod $l^{2}$, for all $k\in $$%
\mathbb{N},$ $k\geq 2.$ Thus, we obtain

\begin{equation*}
A_{n}+A_{n+2^{k}}=\widehat{2}A_{n+2^{k-1}},\text{ for all }k\in \mathbb{N}%
,\;k\geq 2.
\end{equation*}%
ii) Using Proposition 4.9, we have

\begin{equation*}
A_{n}+A_{n+3\cdot 2^{k}}=\widehat{M_{k}\left( M_{k}^{2}-3\right) }%
A_{n+3\cdot 2^{k-1}},k\geq 2.
\end{equation*}%
But $M_{k}\left( M_{k}^{2}-3\right) \equiv 2$ mod $l^{2}$ for all $k\in $$%
\mathbb{N},$ $k\geq 2.$ Therefore, we obtain 
\begin{equation*}
A_{n}+A_{n+3\cdot 2^{k}}=\widehat{2}A_{n+3\cdot 2^{k-1}},\text{ for all }%
k\in \mathbb{N},\;k\geq 2.
\end{equation*}%
iii) 
\begin{equation*}
A_{n}\text{+}A_{n+1}\text{+...+}A_{n+2l^{2}-1}\text{=}%
\sum_{k=0}^{2l^{2}-1}a_{n+k}\text{+}e_{1}\sum_{k=0}^{2l^{2}-1}a_{n+1+k}\text{%
+}e_{2}\sum_{k=0}^{2l^{2}-1}a_{n+2+k}\text{+}e_{3}%
\sum_{k=0}^{2l^{2}-1}a_{n+3+k}.
\end{equation*}%
Since $l$ is odd, applying Proposition 4.4, we obtain: 
\begin{equation*}
\sum_{k=0}^{2l^{2}-1}a_{n+k}=a_{n}+a_{n+1}+...+a_{n+2l^{2}-1}\equiv \underset%
{_{l^{2}-times}}{\underbrace{1+1+...+1}}\text{mod }l\equiv 0~\text{mod }l.
\end{equation*}%
In the same way, we obtain $\overset{2l^{2}-1}{\underset{k=0}{\sum }}%
a_{n+1+k}$$\equiv 0~$mod $l$, $\overset{2l^{2}-1}{\underset{k=0}{\sum }}%
a_{n+2+k}$$\equiv 0~$mod $l$ and $\overset{2l^{2}-1}{\underset{k=0}{\sum }}%
a_{n+3+k}$$\equiv 0~$mod $l$. It results that 
\begin{equation*}
A_{n}+A_{n+1}+...+A_{n+2l^{2}-1}=\widehat{0}~\text{in}~\mathbb{Z}_{l^{2}}.
\end{equation*}%
$\Box \medskip $

\textbf{Conclusions.} In this paper, we gave properties and applications of
some special integer sequences. We generalized Cassini's identity and Pisano
period for a difference equation of degree $k$. Moreover, we provided a new
application in Cryptography and we presented applications of some special
number sequences and quaternion elements over finite rings. The above
results show us that the study of these sequences can provide us new
interesting properties and applications. This remark makes us to continue
their study in further researches. 
\begin{equation*}
\end{equation*}

\textbf{Acknowledgements}. Authors thank organizers of IECMSA-2019 for the
opportunity to present some of their results at this conference.

\begin{equation*}
\end{equation*}

\textbf{References}

\begin{equation*}
\end{equation*}

[EPC; 18] European Payments Council, \textit{Guidelines on Cryptographic
Algorithms Usage and Key Management}, 2018 \qquad

[FP; 09] Falcon, S., Plaza, A., \textit{On k-Fibonacci numbers of arithmetic
indexes}, Applied Mathematics and Computation, 208(2009), 180--185.

\textbf{[}Fib.\textbf{]}
http://www.maths.surrey.ac.uk/hosted-sites/R.Knott/Fibonacci/fib.html

[Fl; 19] C. Flaut, \textit{Some application of difference equations in
Cryptography and Coding Theory}, Journal of Difference Equations and
Applications, 25(7)(2019), 905-920.

\textbf{[}Fl, Sa; 15\textbf{]} C. Flaut, D. Savin, \textit{Quaternion
Algebras and Generalized Fibonacci-Lucas Quaternions}, Adv. Appl. Clifford
Algebras, 25(4)(2015), p. 853-862.

\textbf{[}Fl, Sa; 19\textbf{] }C. Flaut, D. Savin, \textit{Some remarks
regarding l-elements defined in algebras obtained by the Cayley-Dickson
process}, Chaos, Solitons \& Fractals, 118(2019), 112-116.

[FSZ; 19] C. Flaut, D. Savin, G. Zaharia, \textit{Some applications of
Fibonacci and Lucas numbers}, accepted in C. Flaut, S. Hoskova-Mayerova, F.
Maturo, \textit{Algorithms as an approach of applied mathematics}, Springer,
2020.

[Gr, Mi, Ma; 15] J. M. Grau, C. Miguel and A. M. Oller-Marcen, \textit{On
the structure of quaternion rings over} $\mathbb{Z}/n\mathbb{Z},$ Advances
in Applied Clifford Algebras, vol. 25, Issue 4 (2015), p. 875-887.

[Ha; 12] S. Halici, \textit{On Fibonacci Quaternions}, Adv. in Appl.
Clifford Algebras 22(2)(2012), 321-327.

\textbf{[}Ho; 63\textbf{]}\ A. F. Horadam, \textit{Complex Fibonacci Numbers
and Fibonacci Quaternions}, Amer. Math. Monthly, 70(1963), 289--291.

[Jo] R.C. Johnson, \textit{Fibonacci numbers and matrices}, available at%
\newline
http://maths.dur.ac.uk/dma0rcj/PED/fib.pdf.

[Ko; 94] N. Koblitz, \textit{A Course in Number Theory and Cryptography},
Springer Verlag, New-York, 1994, p. 65-76.

[Me; 99] R. Melham, \textit{Sums Involving Fibonacci and Pell Numbers, }%
Portugalie Mathematica, 56(3)(1999), 309-317.

[Re; 13] M. Renault, \textit{The Period, Rank, and Order of the (a,
b)-Fibonacci Sequence Mod m}, Mathematics Magazine, 86(5)(2013), 372-380,%
\newline
https://doi.org/10.4169/math.mag.86.5.372.

[Sa; 17] D. Savin, \textit{About special elements in quaternion algebras
over finite fields}, Advances in Applied Clifford Algebras, vol. 27, June
2017, Issue 2 , p. 1801- 1813.

[Sa; 19] D. Savin, \textit{Special numbers, special quaternions and special
symbol ele- ments}, chapter in the book Models and Theories in Social
Systems, vol. 179, Springer 2019, ISBN-978-3-030-00083-7 , p. 417-430.

[St; 06] Stakhov, A.P., \textit{Fibonacci matrices, a generalization of the
\textquotedblleft Cassini formula\textquotedblright , and a new coding theory%
}, Chaos, Solitons and Fractals, 30(2006), 56-66.

[St; 07] Stakhov, A.P., \textit{The \textquotedblleft
golden\textquotedblright\ matrices and a new kind of cryptography}, Chaos,
Solitons and Fractals, 32(2007), 1138--1146.

[Wa; 60] D. D. Wall, \textit{Fibonacci Series Modulo m}, The American
Mathematical Monthly,67(6)(1960), 525-532.%
\begin{equation*}
\end{equation*}%
\medskip \qquad\ \qquad\ \ 

Cristina FLAUT

{\small Faculty of Mathematics and Computer Science, }\newline

{\small Ovidius University of Constan\c{t}a, Rom\^{a}nia,}

{\small Bd. Mamaia 124, 900527,}

{\small http://www.univ-ovidius.ro/math/}

{\small e-mail: cflaut@univ-ovidius.ro; cristina\_flaut@yahoo.com}%
\begin{equation*}
\end{equation*}%
\medskip \qquad\ \qquad\ \ 

Diana SAVIN

{\small Faculty of Mathematics and Computer Science, }

{\small Ovidius University of Constan\c{t}a, Rom\^{a}nia, }

{\small Bd. Mamaia 124, 900527, }

{\small http://www.univ-ovidius.ro/math/}

{\small e-mail: \ savin.diana@univ-ovidius.ro, \ dianet72@yahoo.com}

\begin{equation*}
\end{equation*}%
\qquad

Geanina ZAHARIA

{\small PhD student at Doctoral School of Mathematics,}

{\small Ovidius University of Constan\c{t}a, Rom\^{a}nia,}

{\small geaninazaharia@yahoo.com}

\end{document}